\documentclass[graybox]{svmult}

\usepackage{type1cm}        % activate if the above 3 fonts are
\usepackage[pdftex]{graphicx}

\usepackage{makeidx}         % allows index generation
\usepackage{graphicx}        % standard LaTeX graphics tool
\usepackage{multicol}        % used for the two-column index
\usepackage[bottom]{footmisc}% places footnotes at page bottom

\usepackage{cite}
\usepackage{newtxtext}       % 
\usepackage{newtxmath}       % selects Times Roman as basic font

\makeindex             % used for the subject index
\newcommand\norm[1]{\left\lVert#1\right\rVert}                    % your makeindex program

\begin{document}
	
	\title*{Coupled Flow and Mechanics in a 3D Porous Media with Line Sources}
	% Use \titlerunning{Short Title} for an abbreviated version of
	% your contribution title if the original one is too long
	\author{Nadia S. Taki and Ingeborg G. Gjerde}
	% Use \authorrunning{Short Title} for an abbreviated version of
	% your contribution title if the original one is too long
	\institute{Nadia S. Taki \at University of Stuttgart, Intitute for Modelling Hydraulic and Environmental Systems, Pfaffenwaldring 5a, 70569 Stuttgart, Germany. \email{nadia.skoglund@iws.uni-stuttgart.de}
		\and Ingeborg G. Gjerde \at University of Bergen, Department of Mathematics,  Natural Science Building, 5007 Bergen, Norway. \email{ingeborg.gjerde@uib.no}}
	%
	% Use the package "url.sty" to avoid
	% problems with special characters
	% used in your e-mail or web address
	%
	\maketitle
	
	\abstract*{In this paper, we consider the numerical approximation of the quasi-static, linear Biot model in a 3D domain $\Omega$ when the right-hand side of the flow equation is concentrated on a 1D line source $\delta_{\Lambda}$. This model is of interest in the context of medicine, where it can be used to model flow and deformation through vascularized tissue. The model itself is challenging to approximate as the line source induces the pressure and flux solution to be singular. To overcome this, we here combine two methods: (i) a fixed-stress splitting scheme to decouple the flow and mechanics equations and (ii) a singularity removal method for the pressure and flux variables. The singularity removal is based on a splitting of the solution into a lower regularity term capturing the solution singularities and a higher regularity term denoted the remainder. With this in hand, the flow equations can now be reformulated so that they are posed with respect to the remainder terms. The reformulated system is then approximated using the fixed-stress splitting scheme. We conclude by showing the results for a test case simulating flow through vascularized tissue. This numerical method is found to converge optimally.}
	
	\abstract{In this paper, we consider the numerical approximation of the quasi-static, linear Biot model in a 3D domain $\Omega$ when the right-hand side of the flow equation is concentrated on a 1D line source $\delta_{\Lambda}$. This model is of interest in the context of medicine, where it can be used to model flow and deformation through vascularized tissue. The model itself is challenging to approximate as the line source induces the pressure and flux solution to be singular. To overcome this, we here combine two methods: (i) a fixed-stress splitting scheme to decouple the flow and mechanics equations and (ii) a singularity removal method for the pressure and flux variables. The singularity removal is based on a splitting of the solution into a lower regularity term capturing the solution singularities and a higher regularity term denoted the remainder. With this in hand, the flow equations can now be reformulated so that they are posed with respect to the remainder terms. The reformulated system is then approximated using the fixed-stress splitting scheme. We conclude by showing the results for a test case simulating flow through vascularized tissue. Here, the  numerical method is found to converge optimally using lowest-order elements for the spatial discretization.}
	
	\section{Introduction}
	\label{sec:1}
	The coupling of mechanics and flow in porous media is relevant for a wide range of applications, occurring for instance in geophysics \cite{bause,erlend} and medicine \cite{simula,Nagashima1987BiomechanicsOH}. We put forward a model relevant for simulating perfusion, i.e., blood flow, and deformation in vascularized tissue. This problem is of high interest in the context of medicine, as clinical measurements of perfusion provide important indicators for e.g. Alzheimer's disease \cite{Chen1977, IturriaMedina2016}, stroke \cite{Markus353} and cancer \cite{Gillies1999}. Moreover, both the tissue and blood vessels are elastic, and these properties constitute another valuable clinical indicator. Vascular compliance, as one example, is reduced in cases of vascular dementia, but not in cases Alzheimer's disease \cite{brain5}. 
	
	We consider the fully coupled quasi-static, linear Biot system \cite{biot1}, modeling a poroelastic media when the source term in the flow equation is concentrated on $\Lambda$. Let $\Omega \subset \mathbb{R}^3$ denote a bounded open 3D domain with smooth boundary $\partial \Omega$ and $\Lambda = \sum_{i=1}^n \Lambda_i$ a collection of straight line segments $\Lambda_i \subset \mathbb{R} \subset \Omega$ embedded in $\Omega$. The equations on the space-time domain $\Omega \times (0, T )$ read:
	
	Find $(\mathbf{u}, p, \mathbf{w})$ such that:
	\begin{subequations}
		\begin{eqnarray}
		-\nabla \cdot [2 \ \mu \ \pmb{\varepsilon}(\mathbf{u})\ +\ \lambda (\nabla \cdot \mathbf{u})\ \mathbf{I}] \ +\ \alpha   \nabla p \ &=& \ \mathbf{f}, \label{eq:h1} \\
		\partial_{t} \bigg(\frac{p}{M} \ +\ \alpha  \nabla \cdot \mathbf{u}\bigg) \ + \ \nabla \cdot \mathbf{w}\  &=&\  \psi + f\delta_{\Lambda}, 	\label{eq:mass} \\
		\kappa^{-1} \ \mathbf{w} \ +\ \nabla p \ &=&\ \rho_{f} \mathbf{g},	\label{eq:h3}
		\end{eqnarray}
	\end{subequations} 
	where $\mathbf{u}$ denotes the displacement, $\pmb{\varepsilon}(\mathbf{u}) = \frac{1}{2}(\nabla \mathbf{u} + \nabla \mathbf{u}^{T})$ the (linear) strain tensor, $p$ the pressure, $\mathbf{w}$ the Darcy's flux, $\alpha$ the Biot coefficient, $\kappa$ the  permeability tensor divided by the fluid viscosity, $\rho_{f}$ the fluid density, $\mathbf{g}$ the gravity vector, $M$ the Biot modulus, $\mu$ and $\lambda$ are Lamé parameters, $\mathbf{f} \in L^2(\Omega)$ the contribution from body forces and a source term $\psi \in L^2(\Omega)$. We assume $M, \alpha, \mu, \lambda \in L^\infty(\Omega)$ to be strictly positive and uniformly bounded and $\rho_{f} \in \mathbb{R}$, $\mathbf{g} \in \mathbb{R}^3$. Additionally, $\kappa \in W^{2,\infty}(\Omega)$ is assumed scalar-valued, as required by the singularity removal method. For simplicity, we use homogeneous boundary conditions $\mathbf{u} = \mathbf{0}$, $p=0$ on $\partial \Omega \times [0, T]$ and initial conditions $\mathbf{u} = \mathbf{u}_{0}$, $p = p_{0}$ in $\Omega \times \{ 0 \}$.  
	
The system \eqref{eq:h1}-\eqref{eq:h3} is made non-standard by a generalized Dirac line source $\delta_{\Lambda}$ of intensity $f$ in the right-hand side of \eqref{eq:mass}. The line source is defined mathematically as
	\begin{equation*}
		\int_{\Omega} f \delta_{\Lambda} v \ \text{d} \Omega = \sum_{j = 1}^{m} \int_{\Lambda_{j}} f(s_{j},t) v(s_{j}) \ \text{d}S \quad \forall v \in C^{0}(\bar{\Omega}).
	\end{equation*}
Physically, it is introduced to model the mass exchange between the vascular network and the surrounding tissue. This exchange occurs through the capillary blood vessels, which have radii ranging from 5 to 10 micrometer. These blood vessels are too small to be captured as 3D objects in a mesh; instead, they are typically reduced to being one-dimensional line segments, see e.g. \cite{1dmodel1, koppl2017, Laurino2019, ingeborg3, vidotto2018}. The system \eqref{eq:h1}-\eqref{eq:h3} would then be on the same form as the one considered in \cite{simula}, with the exception that the exchange term is here concentrated on the 1D vascular network. 

We take $f \in C^0(\Omega)$ and focus our attention on the challenges introduced by the line source $\delta_\Lambda$. The line source induces $p$ and $\mathbf{w}$ to be singular, i.e., they both diverge to infinity on $\Lambda$. Consequently, one has $p \in L^2(\Omega)$ but $\mathbf{w} \notin (L^2(\Omega))^3$; the solution is then not regular enough to fit the analytic framework of \cite{source22,florin}. They both prove global convergence for the fixed-stress splitting scheme applied to Biot's equations. Moreover, these singularities are expensive to resolve numerically, making the solution highly challenging to approximate. 
	
In order to tackle this issue, we here combine two strategies: (i) a fixed-stress splitting scheme that decouples the mechanics equation \eqref{eq:h1} from the flow equations \eqref{eq:mass}-\eqref{eq:h3}, and (ii) a singularity removal method for the flow equations. For an introduction to the fixed-stress splitting scheme, we refer to the works of Mikeli\'c et al. \cite{source22} and Both et al. \cite{florin}; for an introduction to the singularity removal method, we refer to our earlier work \cite{ingeborg, ingeborg2}.

\section{Mathematical Model and Discretization}
In this section, we begin by introducing a splitting method that decomposes $p$ and $\mathbf{w}$ into higher and lower regularity terms. Here, the lower regularity terms are given explicitly. The model \eqref{eq:h1}-\eqref{eq:h3} can be reformulated so that it is given with respect to the higher regularity terms; we refer to this as the singularity removal method. Next, we show how this model can readily be approximated by the fixed-stress splitting scheme. 
	\iffalse
	Let $d\in \{2,3\}$ denote the dimension and $ W^{2,\infty} (\Omega) = \{  u\in L^{\infty}(\Omega) :  D^{\beta}u \in L^{\infty}(\Omega) \text{ for } |\beta| \leq 2  \}$ with $D^{\beta}$ as the weak distributional derivative of $u$ with multiindex $\beta$. We then make these assumptions used on the effective coefficients:
	\begin{enumerate}
		\item Let $\rho_{f} \in \mathbb{R}$, $\mathbf{g} \in \mathbb{R}^{3}$ be constant.
		\item Let $M, \alpha, \mu, \lambda \in L^{\infty} (\Omega)$ be positive, uniformly bounded, with the lower bounded strictly positive.
		\item Let $\kappa \in W^{2,\infty} (\Omega)^{d\times d}$ be a symmetric matrix, which is constant in time and has uniformly bounded eigenvalues, i.e., there exist constants $k_{m},k_{M} \in \mathbb{R}$, satisfying for all $\pmb{x} \in \Omega$ and for all $\mathbf{z} \in \mathbb{R}^{3} \backslash \{0\}$. 
		\begin{equation*}
		0 < k_{m}\ \mathbf{z}^{T}\mathbf{z} \leq \mathbf{z}^{T} \kappa(x) \ \mathbf{z} \leq k_{M} \ \mathbf{z}^{T}\mathbf{z} < \infty.
		\end{equation*}
	\end{enumerate}
	\fi
	
	\subsection{Singularity Removal Method}\label{singularremoval}
	%In this section, we will show how a solution splitting can be used to reformulate the flow equations \eqref{eq:mass}-\eqref{eq:h3} so that they are given with respect to higher regularity terms. 
	For the sake of notational simplicity, we assume $\kappa$ to be constant; a spatially varying $\kappa$ could be handled as shown in \cite[Sect. 3.3]{ingeborg}. Let $\mathbf{a}_i, \mathbf{b}_i$ denote the endpoints of the line segment $\Lambda_i$. From \cite[Sect. 3.2]{ingeborg}, we have a function $G$ defined as 
	\begin{equation}\label{eq:G}
	G(\pmb{x}) = \sum_{i=1}^n \frac{1}{4 \pi} \ln\bigg (\frac{r_{b,i} + L_i + \pmb{\gamma_i} \cdot (\pmb{a}_i - \pmb{x})}{r_{a,i}  + \pmb{\gamma}_i \cdot (\pmb{a}_i - \pmb{x})} \bigg),
	\end{equation}
	with $r_{a,i} = \norm{\pmb{x}- \pmb{a}_i}$, $r_{b,i} = \norm{\pmb{x}- \pmb{b}_i}$, $L_i = \norm{\pmb{b}_i - \pmb{a}_i}$ and $\pmb{\gamma}_i = \frac{\pmb{b}_i-\pmb{a}_i}{L_i}$ as the normalized tangent vector of $\Lambda_i$. Centrally, this function solves $-\Delta G = \delta_\Lambda$ in the appropriate weak sense; i.e., we have:
	\begin{equation*} 
	-\int_{\Omega} \Delta G \ v \ \text{d} \Omega = \int_{\Lambda} \ v \ \text{d}S \quad \forall v \in C^{0}(\bar{\Omega}).
	\end{equation*}
	Having this function in hand, we next formulate the following splitting ansatz:
	\begin{equation}
	p = p_{s} + p_{r}, \quad \mathbf{w} = \mathbf{w}_{s} + \mathbf{w}_{r},\label{eq:phelp}
	\end{equation}
	where $p_{s} = \frac{f G(\pmb{x})}{\kappa}$ and $\mathbf{w}_{s} = - \kappa\nabla p_{s}$. The terms $p_s$ and $\mathbf{w}_{s}$ capture the singular part of the solution, and are explicitly given via the function $G$. This allows $p_r$ and $\mathbf{w}_r$ to enjoy higher regularity and improved approximation properties. Assume for the moment that the solution $\mathbf{u}$ is given. Inserting the splitting \eqref{eq:phelp} into \eqref{eq:mass}-\eqref{eq:h3} one finds the following reformulated flow equation:
	
	Find $(p_r, \mathbf{w}_r)$ such that:
	\begin{subequations}
		\begin{eqnarray}
		\partial_{t} \bigg(\frac{p_r}{M} \ +\ \alpha  \nabla \cdot \mathbf{u}\bigg) \ + \ \nabla \cdot \mathbf{w}_r \  &=&\  \psi_{r}, 	\label{eq:mass-reform} \\
		\kappa^{-1} \ \mathbf{w}_r \ +\ \nabla p_r \ &=&\ \rho_{f} \mathbf{g},	\label{eq:h3-reform}
		\end{eqnarray}
	\end{subequations}
	where $\psi_{r} = \psi - \frac{\partial_t p_s}{M} +  G \Delta f + 2 \nabla G \cdot \nabla f$.
	Here, \eqref{eq:h3-reform} is straightforward to obtain. For \eqref{eq:mass-reform}, we used that 
	\begin{eqnarray*}
	\partial_{t} \bigg(\frac{p_r}{M} \ +\ \alpha  \nabla \cdot \mathbf{u}\bigg) \ +  \nabla \cdot \mathbf{w}_r \  &=& \psi + f\delta_\Lambda - \frac{\partial_t p_s}{M} - \nabla \cdot \mathbf{w}_s \\
	&=& \psi + f\delta_\Lambda - \frac{\partial_t p_s}{M} + \nabla \cdot (\kappa \nabla \frac{f G}{\kappa}) \\
	&=& \psi - \frac{\partial_t p_s}{M} + 2\nabla f \cdot \nabla G + (\Delta f) G.
	\end{eqnarray*}
	In the last line we used the product rule to obtain $\nabla \cdot (\kappa \nabla \frac{fG}{\kappa})= \Delta (fG)=f \Delta G + 2\nabla f \cdot \nabla G+ (\Delta f) G$ along with the relation $f\Delta G = -f \delta_\Lambda$.
	
The value of the reformulation lies in the fact that $\psi_r$ can now be expected to belong to $L^2(\Omega)$. To see this, note that $\psi \in L^2(\Omega)$ by assumption. $G \in L^2(\Omega)$ can be shown by straightforward calculation; it follows that $p_s \in L^2(\Omega)$. Finally, one can show that $\nabla G \cdot \nabla f \in L^2(\Omega)$; for verification of this, we refer to the calculations in \cite[Sect. 4.2]{ingeborg2} along with the embedding $f \in C^0(\Omega) \subset H^1(\Omega)$. 

Let now \eqref{eq:h1} and \eqref{eq:mass-reform}-\eqref{eq:h3-reform} denote the reformulated Biot equation. As $\psi_r \in L^2(\Omega)$, this system fits the analytic framework of \cite{florin}.
	\subsection{Fixed-Stress Splitting Scheme}
	\label{singremfsss}
	Next, we show how the reformulated Biot equations \eqref{eq:h1} and \eqref{eq:mass-reform}-\eqref{eq:h3-reform} can be approximated via the fixed-stress splitting scheme from \cite{florin}.
	Let $\mathcal{T}_{h}$ be the triangularization of the domain $\Omega$ with mesh size $h$. We let $0=t^{0} < \cdots < t^{N} = T$ be a partition of the time interval $(0,T)$ with $N \in \mathbb{N}^{*}$ and define a constant time step size $\tau = t^{k+1} - t^{k} := T / N$ for $k \geq 0$. To discretize the system, we employ backward Euler for time and a finite element method for space. The solutions are approximated with linear piecewise polynomials, constant piecewise polynomials and
	lowest-order Raviart-Thomas spaces for the displacement, pressure and flux,
	respectively. The discrete spaces are given by:
	\begin{eqnarray*}
		\mathbf{V}_{h} &=& \{ \ \mathbf{v}_{h} \in [ \ H^{1}_{0}(\Omega)\ ]^{3} \ | \ \forall K  \in  \mathcal{T}_{h}, \ \mathbf{v}_{h\vert K} \in [\ \mathbb{P}_{1} \ ]^{3} \}, \\ 
		Q_{h} &=& \{ \ q_{h} \in L^{2}(\Omega) \ |\ \forall K \in  	\mathcal{T}_{h}, \ q_{h\vert K} \in \mathbb{P}_{0} \}, \\
		\mathbf{Z}_{h} &=& \{ \mathbf{z}_{h} \in H(\mathrm{div};\Omega)\ |\ \forall 	K \in  \mathcal{T}_{h}, \ \mathbf{z}_{h\vert K}(\pmb{x}) = \pmb{\eta}+\xi\pmb{x}, \ \pmb{\eta} \in \mathbb{R}^{3},\xi \in \mathbb{R} \},
	\end{eqnarray*}
	with $\mathbb{P}_{1}$ and $\mathbb{P}_{0}$ as the linear and constant piecewise polynomials.
	
	Take now $\langle \cdot, \cdot \rangle$ to be the $L^{2}(\Omega)$-inner product and $(\mathbf{u}^{0}_{h}, p^{0}_{h}, \mathbf{w}^{0}_{h}) \in \mathbf{V}_{h} \times Q_{h} \times \mathbf{Z}_{h}$ to be the initial values of the solution. We assume the solution of the displacement, pressure and flux is known for the previous time step. The time-discretization of \eqref{eq:h1} and \eqref{eq:mass-reform}-\eqref{eq:h3-reform} then reads:
	
	Given $(\mathbf{u}^{n-1}_{h}, p^{n-1}_{r,h},\mathbf{w}^{n-1}_{r, h}) \in \mathbf{V}_{h} \times Q_{h} \times \mathbf{Z}_{h}$, find $(\mathbf{u}^{n}_{h}, p^{n}_{r,h}, \mathbf{w}^{n}_{r,h}) \in \mathbf{V}_{h} \times Q_{h} \times \mathbf{Z}_{h}$ such that 
	\begin{subequations}
		\begin{eqnarray*}
			\langle 2 \mu \pmb{\varepsilon}(\mathbf{u}^{n}_{h}), \pmb{\varepsilon}(\mathbf{v}_{h})  \rangle + \langle \lambda (\nabla \cdot \mathbf{u}^{n}_{h}), \nabla \cdot \mathbf{v}_{h} \rangle - \langle \alpha p^{n}_{h}, \nabla \cdot \mathbf{v}_{h} \rangle &=& \langle \mathbf{f}^{n}, \mathbf{v}_{h} \rangle,  \\
			\bigg \langle \frac{1}{M} p^{n}_{r,h}, q_{h} \bigg\rangle + \langle \alpha \nabla \cdot \mathbf{u}^{n}_{h}, q_{h} \rangle + \tau \langle \nabla \cdot \mathbf{w}^{n}_{r,h}, q_{h} \rangle  &=& \tau \langle \psi^{n}_r, q_{h} \rangle+ \bigg \langle \frac{1}{M} p^{n-1}_{r,h}, q_{h} \bigg\rangle \\
			&  +& \langle \alpha \nabla \cdot \mathbf{u}^{n-1}_{h}, q_{h} \rangle,
			\\ 
			\langle \kappa^{-1} \mathbf{w}^{\ n}_{r,h}, \mathbf{z}_{h} \rangle - \langle p^{n}_{r,h}, \nabla \cdot \mathbf{z}_{h} \rangle &=& \langle \rho_{f} \mathbf{g}, \mathbf{z}_{h} \rangle,
		\end{eqnarray*}
	\end{subequations}
	for all $(\mathbf{v}_{h}, q_{h}, \mathbf{z}_{h}) \in \mathbf{V}_{h} \times Q_{h} \times \mathbf{Z}_{h}$. 
	
	The idea of the fixed-stress splitting scheme is to decouple the flow and mechanics equation while keeping an artificial volumetric stress $\sigma_{\beta} = \sigma_{0} + K_{dr} \nabla \cdot \mathbf{u} - \alpha p$ constant. Here, $K_{dr} \in L^{\infty}(\Omega)$ is referred to as the drained bulk modulus. We consider the theoretically optimal tuning parameter $\beta_{FS} = \alpha^{2}/K_{dr}$ with $K_{dr} = \frac{d}{2}(\mu + \lambda)$ \cite{florin}.\\
	\indent
	We define a sequence $(\mathbf{u}^{n,i}_{h}, \ p^{n,i}_{r,h}, \  \mathbf{w}^{n,i}_{r,h})$, $i \geq 0$. Let $i$ denote the current iteration step and $i-1$ denote the previous iteration step. Then initialize $\mathbf{u}$, $p_r$ and $\mathbf{w}_r$ by $\mathbf{u}^{n,0}_{h}= \mathbf{u}^{n-1}_{h}$, $p^{n,0}_{r,h}= p^{n-1}_{r,h}$ and $\mathbf{w}^{n,0}_{r,h} = \mathbf{w}^{n-1}_{r,h}$, respectively. The algorithm iterates until a stopping criterion is reached. The full scheme reads:\\
	\\
	\textbf{Step 1:} Given $(\mathbf{u}^{n,i-1}_{h}, \ p^{n,i-1}_{r,h}, \ \mathbf{w}^{n,i-1}_{r,h}) \in \mathbf{V}_{h} \times Q_{h} \times \mathbf{Z}_{h}$.  Find $(p^{n,i}_{h}, \ \mathbf{w}^{n,i}_{h}) \in Q_{h} \times \mathbf{Z}_{h}$ such that $\forall \, (q_{h}, \ \mathbf{z}_{h}) \in Q_{h} \times \mathbf{Z}_{h}$:
	\begin{subequations}
		\begin{eqnarray}\label{eq:FSSSs2}
		\bigg \langle \bigg( \frac{1}{M} + \beta_{FS} \bigg) p^{n,i}_{r,h}, q_{h} \bigg \rangle + \tau \langle \nabla \cdot \mathbf{w}_{r,h}^{n,i}, q_{h} \rangle & =&  \tau \langle \psi^{n}_r, q_{h} \rangle  + \bigg \langle \frac{1}{M} p^{n-1}_{r,h}, q_{h} \bigg \rangle \\
		& +& \langle \alpha \nabla \cdot \mathbf{u}^{n-1}_{h}, q_{h} \rangle  + \langle \beta_{FS} \ p^{n,i-1}_{r,h}, q_h \rangle \notag\\
		&-& \langle \alpha \nabla \cdot \mathbf{u}^{n,i-1}_{h}, q_{h} \rangle,\notag\\
		\langle \kappa^{-1} \mathbf{w}_{r,h}^{n,i}, \mathbf{z}_{h} \rangle - \langle p^{n,i}_{r,h}, \nabla \cdot \mathbf{z}_{h} \rangle &=&  \langle \rho_{f} \mathbf{g},\label{eq:FSSSs1} \mathbf{z}_{h} \rangle.
		\end{eqnarray}
		\textbf{Step 2:} Update the full pressure and flux solutions: $p^{n,i}_{h} = p^{n}_{s,h} + p^{n,i}_{r,h}$ and $\mathbf{w}^{n,i}_{h} = \mathbf{w}^{n}_{s,h} + \mathbf{w}^{n,i}_{r,h}$. \\
		\textbf{Step 3:} Given  $p^{n,i}_{h} \in Q_{h}$. Find $\mathbf{u}_{h}^{n,i} \in \mathbf{V}_{h}$, such that $\forall \, \mathbf{v}_{h} \in \mathbf{V}_{h}$:
		\begin{equation}\label{eq:FSSS2}
		\langle 2 \mu \pmb{\varepsilon}(\mathbf{u}^{n,i}_{h}), \pmb{\varepsilon}(\mathbf{v}_{h})  \rangle + \langle \lambda (\nabla \cdot \mathbf{u}^{n,i}_{h}), \nabla \cdot \mathbf{v}_{h} \rangle = \langle \alpha p^{n,i}_{h}, \nabla \cdot \mathbf{v}_{h} \rangle + \langle \mathbf{f}^{n}, \mathbf{v}_{h} \rangle.
		\end{equation}
	\end{subequations} 
	\section{Numerical Results}\label{numres}
	In this section, we provide numerical convergence results for a test case using parameters relevant for flow through vascularized tissue. 
	Let the medium of consideration be an isotropic, homogeneous porous medium and $\kappa$ a positive scalar quantity. We let $(\mathbf{p}_{r,h}^{i} , \mathbf{w}_{r,h}^{i}, \mathbf{u}_{h}^{i})$ be the solutions at iteration step $i$ and $(\mathbf{p}_{r,h}^{i-1} , \ \mathbf{w}_{r,h}^{i-1}, \ \mathbf{u}_{h}^{i-1})$ the solutions at the previous iteration step $i-1$. The procedure stops when reaching the following criterion:
	\begin{equation*}
		\norm{(\mathbf{p}_{r,h}^{i} , \ \mathbf{w}_{r,h}^{i}, \ \mathbf{u}_{h}^{i}) - (\mathbf{p}_{r,h}^{i-1} , \ \mathbf{w}_{r,h}^{i-1}, \ \mathbf{u}_{h}^{i-1})} \ \leq \ \epsilon_{a} + \epsilon_{r} \  \norm{(\mathbf{p}_{r,h}^{i} , \ \mathbf{w}_{r,h}^{i}, \ \mathbf{u}_{h}^{i})},
	\end{equation*} 
	where $\epsilon_{a}, \epsilon_{r} > 0$ are given tolerances (see Table \ref{table:valid2}). 
	\begin{table}[t]
		\centering
		\caption{Material parameters used to solve the Biot's equations with lower dimensional source terms \eqref{eq:h1}-\eqref{eq:h3} in Sect. \ref{numres}. There is a wide rage of parameters used in literature. The ones represented here are a sample of representative parameters.}
		\begin{tabular}{@{}llll@{}}
			\hline
			Symbol: & Quantity: & Value: & Reference:\\
			\hline
			$\kappa$ & Permeability divided&  $1.57e-2$ mm$^2$ mPa$^{-1}$ s$^{-1}$ &\cite[Table 1]{Weiss2017}\\
			&  by the fluid viscosity & &\\
			$E$  & Tuning parameter & $1.5e6$ mPa& \cite[Table 6]{simula}\\
			$M$ & Biot modulus & $3.9e7$ mPa &\cite[Table 2]{Guo2017} \\
			$\alpha$ & Biot coefficient & $1.0$ &\\
			$\nu$ & Poisson's ratio & $0.2$ & \\
			$\mathbf{g}$ & Gravitational vector & $\mathbf{0}$ mm s$^{-2}$& \cite{biotbook}\\
			$T$ & Final time & $1.0$ s &\\
			$\tau$  & Time step & $0.1$ s &\\
			$\epsilon_{a}$  & Absolute error tolerance & $1e-6$ &\\
			$\epsilon_{r}$  & Relative error tolerance  & $1e-6$ &\\
			\hline
		\end{tabular}
		\label{table:valid2}
	\end{table}
	
	Let $\Omega = \{(0,1)\times (0,1) \times (0,1) \subset \mathbb{R}^{3} \}$ be a cube discretized by $1/h \times 1/h \times 1/h$ tetrahedrons. The numerical results are obtained by the fixed-stress splitting scheme proposed in Sect. \ref{singremfsss} and the programming platform FEniCS \cite{LoggMardalEtAl2012a}. Convergence is tested against the following analytic solutions
	\begin{equation*}
	p_{r} = \frac{1}{4 \pi \kappa} f(t)( r_{a} - r_{b} ), \quad \quad
	\mathbf{w}_{r} = -\kappa  \nabla p_{r},  \quad \quad
	\mathbf{u} =t x(1-x)y(1-y)z(1-z)  [1\ \ 1\ \ 1]^{T},
	\end{equation*}
	where $f(t) = \sin (t)$ is a pulsative intensity function. We selected two points $\pmb{a} = [0.5 \ \ 0.8 \ \  0.5]^{T}$ and $\pmb{b} = [0.5\ \ 0.2 \ \ 0.5]^{T}$ to describe the line segment. Then computed the solutions using mixed-finite element formulations for the correction terms $p_{r}$ and $\mathbf{w}_{r}$, and the solution displacement $\mathbf{u}$ is calculated with the conformal finite element formulation.
	
\begin{figure}[t]
		\sidecaption
		\includegraphics[scale=.2]{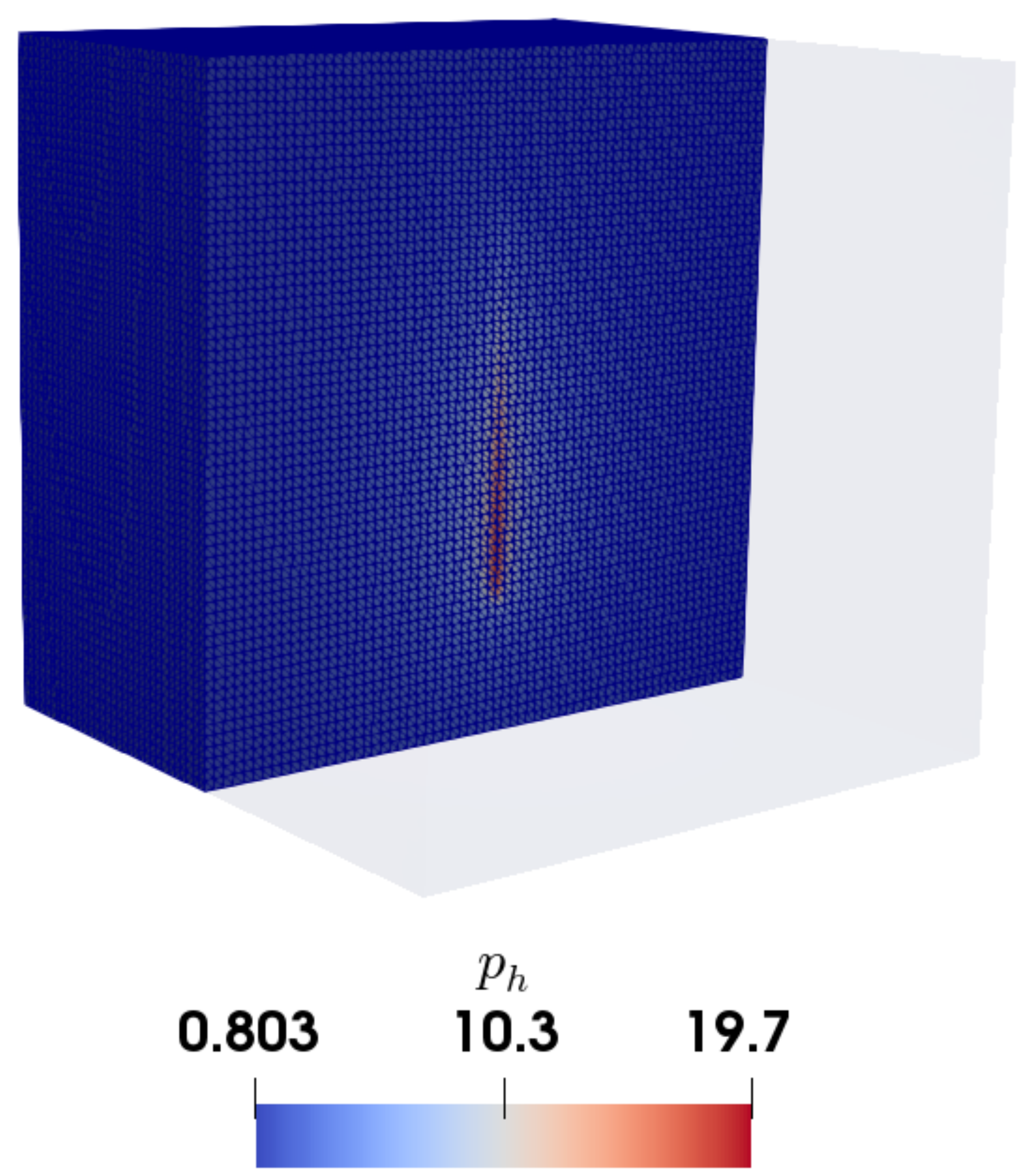}
		\hspace{-6pt}
		\includegraphics[scale=.194]{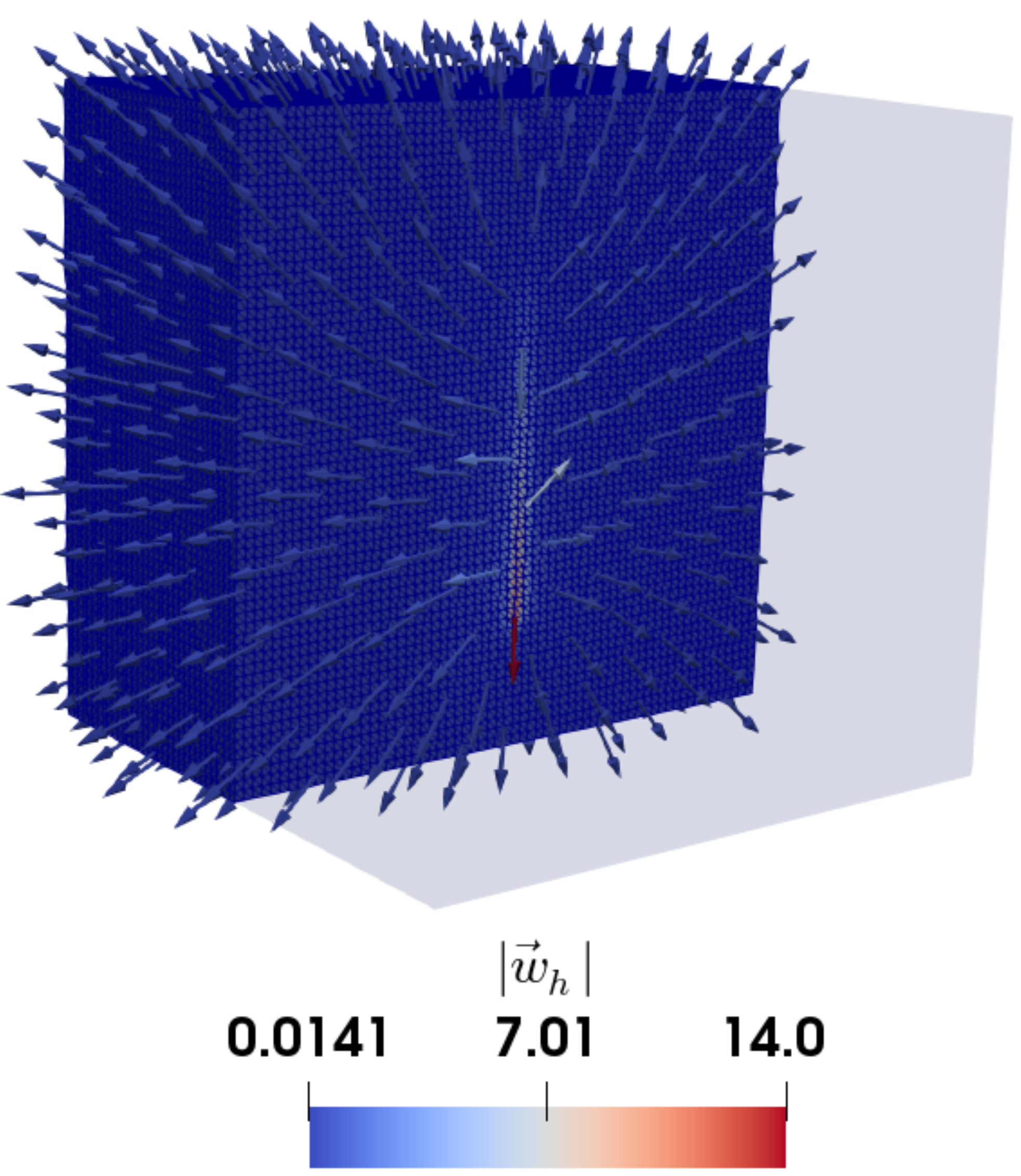}
		\hspace{-6pt}
		\includegraphics[scale=.21]{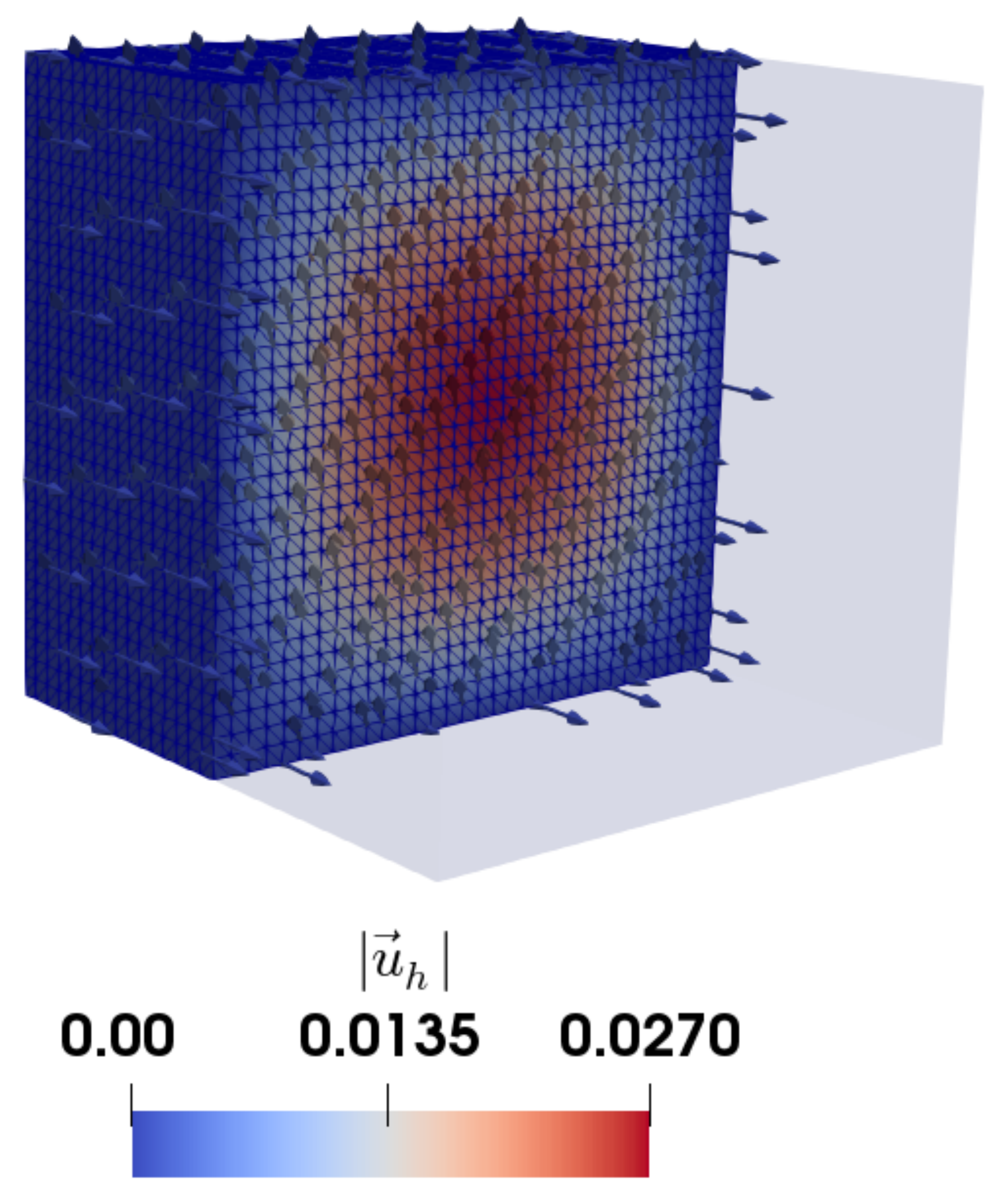}
		\caption{Left: Plot of the reconstructed pressure. Middle: Magnitude of the full flux. Right: Magnitude of the displacement. All plots are numerical solutions obtained by the fixed-stress splitting scheme \eqref{eq:FSSSs2}-\eqref{eq:FSSS2} with one line source.}
		\label{fig:test3d1}
	\end{figure}

	Table \ref{table:valid} shows the error and convergence rates obtained using the parameters listed in Table \ref{table:valid2}. The singularity removal based fixed-stress splitting scheme is seen to converge optimally for each variable $\mathbf{u}_h$, $p_h$, and $\mathbf{w}_h$.
	The plots for this problem are illustrated in Fig. \ref{fig:test3d1}. The figure includes the plot of the full pressure, the magnitude of the full flux and the magnitude of the displacement, accordingly.
	\begin{table}[!t]
		\centering
		\caption{Errors and convergence rates obtained solving \eqref{eq:FSSSs2}-\eqref{eq:FSSS2} with analytical solutions found in this section. For reference, optimal convergence rates are listed in the bottom row.}
		\begin{tabular}{lccc}
			\hline
			{$h$} &      $\norm{{p}_{a} - {p}_{h}}_{L^{2}(\Omega)}$ &  $\norm{\mathbf{w}_{a} - \mathbf{w}_{h}}_{L^{2}(\Omega)}$ &   $\norm{\mathbf{u}_{a} - \mathbf{u}_{h}}_{L^{2}(\Omega)}$ \\
			\hline
			$1/8$     & $1.2e-01$  & $7.2e-03$  & $5.9e-04$  \\
			$1/16$     & $6.3e-02$  & $3.5e-03$  & $1.5e-04$ \\
			$1/32$     & $3.1e-02$  & $1.7e-03$ & $3.7e-05$  \\
			\hline
			Rate &    $1.0$ &   $1.0$ &$ 2.0$ \\
			\hline
			Optimal &    $1.0$ &   $ 1.0$ &$ 2.0$ \\
			\hline
		\end{tabular} \label{table:valid}
	\end{table}
	\bibliographystyle{spmpsci}
	\bibliography{references}

\begin{thebibliography}{10}
\providecommand{\url}[1]{{#1}}
\providecommand{\urlprefix}{URL }
\expandafter\ifx\csname urlstyle\endcsname\relax
  \providecommand{\doi}[1]{DOI~\discretionary{}{}{}#1}\else
  \providecommand{\doi}{DOI~\discretionary{}{}{}\begingroup
  \urlstyle{rm}\Url}\fi

\bibitem{bause}
Bause, M., Radu, F.A., K\"ocher, U.: Space-time finite element approximation of
  the {B}iot poroelasticity system with iterative coupling.
\newblock Computer Methods in Applied Mechanics and Engineering \textbf{320},
  745 -- 768 (2017)

\bibitem{biot1}
Biot, M.A.: General theory of three‐dimensional consolidation.
\newblock Journal of Applied Physics \textbf{12}(2), 155 -- 164 (1941)

\bibitem{florin}
Both, J.W., Borregales, M., Kumar, K., Nordbotten, J.M., Radu, F.A.: Robust
  fixed stress splitting for {B}iot's equations in heterogeneous media.
\newblock Applied Mathematics Letters \textbf{68}, 101 -- 108 (2017)

\bibitem{Chen1977}
Chen, Y., Wolk, D., Reddin, J., Korczykowski, M., Martinez, P., Musiek, E.,
  Newberg, A., Julin, P., Arnold, S., Greenberg, J., Detre, J.: Voxel-level
  comparison of arterial spin-labeled perfusion mri and fdg-pet in alzheimer
  disease.
\newblock Neurology \textbf{77}(22), 1977--1985 (2011)

\bibitem{1dmodel1}
D'Angelo, C., Quarteroni, A.: On the coupling of 1{D} and 3{D}
  diffusion-reaction equations: Application to tissue pefusion problems.
\newblock Mathematical Models and Methods in Applied Sciences \textbf{18}(8),
  1481--1504 (2008)

\bibitem{brain5}
Dhoat, S., Ali, K., Bulpitt, C.J., Rajkumar, C.: Vascular compliance is reduced
  in vascular dementia and not in alzheimer's disease.
\newblock Age and Ageing \textbf{37}(6), 653--659 (2008)

\bibitem{biotbook}
Formaggia, L., Quarteroni, A., Veneziani, A.: Cardiovascular Mathematics:
  Modeling and Simulation of the Circulatory System, vol.~1 (2009)

\bibitem{Gillies1999}
Gillies, R.J., Schomack, P.A., Secomb, T.W., Raghunand, N.: Causes and effects
  of heterogeneous perfusion in tumors.
\newblock Neoplasia \textbf{1}(3), 197--207 (1999)

\bibitem{ingeborg3}
Gjerde, I., Kumar, K., Nordbotten, J.M.: A singularity removal method for
  coupled 1d-3d flow models (2018).
\newblock ArXiv:1812.03055 [math.AP]

\bibitem{ingeborg2}
Gjerde, I., Kumar, K., Nordbotten, J.M.: A mixed approach to the poisson
  problem with line sources (2019).
\newblock ArXiv:1910.11785 [math.AP]

\bibitem{ingeborg}
Gjerde, I., Kumar, K., Nordbotten, J.M., Wohlmuth, B.: Splitting method for
  elliptic equations with line sources.
\newblock ESAIM: M2AN \textbf{53}(5) (2019)

\bibitem{Guo2017}
Guo, L., Vardakis, J., Lassila, T., Mitolo, M., Ravikumar, N., Chou, D., Lange,
  M., Sarrami-Foroushani, A., Tully, B., Taylor, Z., Varma, S., Venneri, A.,
  Frangi, A., Ventikos, Y.: Subject-specific multiporoelastic model for
  exploring the risk factors associated with the early stages of alzheimer's
  disease.
\newblock Interface Focus \textbf{8}(1) (2017)

\bibitem{IturriaMedina2016}
Iturria-Medina, Y., , Sotero, R.C., Toussaint, P.J., Mateos-P{\'{e}}rez, J.M.,
  Evans, A.C.: Early role of vascular dysregulation on late-onset alzheimer's
  disease based on multifactorial data-driven analysis.
\newblock Nature Communications \textbf{7}(1) (2016)

\bibitem{koppl2017}
Köppl, T., Vidotto, E., Wohlmuth, B., Zunino, P.: Mathematical modeling,
  analysis and numerical approximation of second-order elliptic problems with
  inclusions.
\newblock Mathematical Models and Methods in Applied Sciences \textbf{28}(05),
  953--978 (2018)

\bibitem{Laurino2019}
Laurino, F., Zunino, P.: Derivation and analysis of coupled {PDEs} on manifolds
  with high dimensionality gap arising from topological model reduction.
\newblock {ESAIM}: Mathematical Modelling and Numerical Analysis
  \textbf{53}(6), 2047--2080 (2019)

\bibitem{simula}
Lee, J., Piersanti, E., Mardal, K.A., Rognes, M.: A mixed finite element method
  for nearly incompressible multiple-network poroelasticity.
\newblock SIAM Journal on Scientific Computing \textbf{41}(2), A722--A747
  (2019)

\bibitem{LoggMardalEtAl2012a}
Logg, A., Mardal, K.A., Wells, G.N., et~al.: Automated Solution of Differential
  Equations by the Finite Element Method.
\newblock Springer (2012)

\bibitem{Markus353}
Markus, H.S.: Cerebral perfusion and stroke.
\newblock Journal of Neurology, Neurosurgery \& Psychiatry \textbf{75}(3),
  353--361 (2004)

\bibitem{source22}
Mikeli{\'{c}}, A., Wheeler, M.F.: Convergence of iterative coupling for coupled
  flow and geomechanics.
\newblock Computational Geosciences \textbf{17}(3), 455--461 (2013)

\bibitem{Nagashima1987BiomechanicsOH}
Nagashima, T., Tamaki, N., Matsumoto, S., Horwitz, B., Seguchi, Y.:
  Biomechanics of hydrocephalus: a new theoretical model.
\newblock Neurosurgery \textbf{21}(6), 898--904 (1987)

\bibitem{erlend}
Storvik, E., Both, J.W., Kumar, K., Nordbotten, J.M., Radu, F.A.: On the
  optimization of the fixed-stress splitting for {B}iot's equations.
\newblock International Journal for Numerical Methods in Engineering
  \textbf{120}(2), 179--194 (2019)

\bibitem{vidotto2018}
Vidotto, E., Koch, T., K\"{o}ppl, T., Helmig, R., Wohlmuth, B.: Hybrid models
  for simulating blood flow in microvascular networks.
\newblock Multiscale Modeling {\&} Simulation \textbf{17}(3), 1076--1102 (2019)

\bibitem{Weiss2017}
Weiss, C.: Finite element analysis for model parameters distributed on a
  hierarchy of geometric simplices.
\newblock GEOPHYSICS \textbf{82}(4), 1--52 (2017)

\end{thebibliography}
	
\end{document}